\documentclass[12pt,a4paper,reqno]{amsart}
\usepackage{amsmath}
\usepackage{amsfonts}
\usepackage{amssymb}
\numberwithin{equation}{section}

\theoremstyle{plain}

\theoremstyle{definition}

\renewcommand{\leq}{\leqslant}
\renewcommand{\geq}{\geqslant}

\newsavebox{\proofbox}
\savebox{\proofbox}{\begin{picture}(7,7)%
  \put(0,0){\framebox(7,7){}}\end{picture}}

\def\ni{\noindent}

\def\emph#1{{\it #1}}
\def\textbf#1{{\bf #1}}

\begin{document}

\title{On the Lacunarity of some eta-products}

\author{Yudong Wang, Chunlei Liu, Haobo Dai}
\address{Department of Mathematics\\
Shanghai Jiao Tong University\\
Shanghai 200240\\}

\email{wyd007001@sjtu.edu.cn}

\thanks{This work is supported by Project 11071160 of the Natural Science Foundation of China.}


\begin{abstract}
\ni The lacunarity is an interesting property of a formal series. We say a series is lacunary if ''almost all'' of its coefficients are zero. In this article we considered about the lacunarity of some eta-products like $\eta(z)^2\eta(bz)^2$, and proved that they are lacunary if and only if $b\in\{1,2,3,4,16\}$. Then We write them as linear combinations of some CM forms.
\end{abstract}

\maketitle

\section{Introduction}

The Dedekind eta-function, defined as:
\begin{displaymath}
\eta(z):=q^{\frac{1}{24}}\prod_{n=1}^{\infty}(1-q^n)
\end{displaymath}
where $q=2\pi i z$, plays an important role in modular form theories. Moreover, as the generating function of the Partition Function $p(n)$:
\begin{displaymath}
\sum_{n=0}^{\infty}p(n)q^n=\prod_{n=1}^{\infty}(1-q^n)^{-1}
\end{displaymath}
$\eta(z)$ is also very useful in partition theory and combinatorial theory. So are the eta-quotients, defined as the products and quotients of some eta functions:
\begin{displaymath}
f(z):=\prod_{\delta}\eta(\delta z)^{r_{\delta}}=q^s\prod_{n=1}^{\infty}\prod_{\delta}(1-q^{\delta n})^{r_{\delta}}
\end{displaymath}
where $\delta$ are positive integers, and $r_{\delta}$ ,related to $\delta$, are non-zero integers. $s:=\frac{1}{24}\sum_{\delta}\delta r_{\delta}$. When all $r_{\delta}$ are positive, $f(z)$ is also called eta-products.

An interesting property of a series $\sum_{n=0}^{\infty}a(n)q^n$ is the "lacunarity". A series $\sum_{n=0}^{\infty}a(n)q^n$ is said to be $lacunary$ if "almost all" of its
coefficients are zero, that is:
\begin{displaymath}
\lim_{X\rightarrow\infty}\frac{\sharp\{0\leq n \leq X: a(n)=0\}}{X}=1
\end{displaymath}
By far, there have been many researches on the lacunarity of eta-quotients. For example
in [1], Serre proved that given a positive even integer r, $\eta(z)^r$ is lacunary if and only if $r\in\{2,4,6,8,10,14,26\}$. In [2], Clader discussed the eta-quotients in the
form $\frac{\eta(az)^b}{\eta(z)}$, and proved the only 19 cases when they are lacunary. In [3], Martin researched the lacunarity of eta-quotients which are Hecke eigenforms. In this
article, we discuss the eta-products in the form $\eta(z)^2\eta(kz)^2$, and prove the following theorem:

\textbf{Theorem 1}  Let $f_b(z):=\eta(z)^2\eta(bz)^2$, $b$ is a positive integer and has no square integer divisor except 2,3. $23\nmid b$. Then $f_b(z)$ is lacunary if and only if $b\in\{1,2,3,4,16\}$.

We will prove this theorem in section 3, using the theory of modular form and CM forms. After that, we will give the method of computing the coefficients of $f_b(z)$ when it is lacunary
 in section 4.
~\\

\section{Preliminaries}
~\\

The Dedekind eta-function $\eta(24z)$, as we know, is a weight $\frac{1}{2}$ modular form on $\Gamma_0(576)$, with Nebentypus character $\chi_{12}(n):=\left(\frac{12}{n}\right)$  . As consequence, the eta-quotient may be modular form on some modular group.
In fact, we have the following theorem from Ono[4]:

\textbf{Theorem 2.1} If $f(z)=\prod_{\delta |N}\eta(\delta z)^{r_{\delta}}$ is an eta-quotient with $k=\frac{1}{2}\sum_{\delta | N}r_{\delta}$, with the additional properties that
\begin{displaymath}
\sum_{\delta | N}\delta r_{\delta}\equiv 0~(mod~24)
\end{displaymath}
and
\begin{displaymath}
\sum_{\delta | N}\frac{N}{\delta} r_{\delta}\equiv 0~(mod~24)
\end{displaymath}
then $f(z)$ satisfies
\begin{displaymath}
f\left(\frac{az+b}{cz+d}\right)=\chi(d)(cz+d)^k f(z)
\end{displaymath}
for every $\left(\begin{array}{ccc} a&b\\c&d \end{array}\right)\in \Gamma_0(N)$. Here the character $\chi$ is defined by $\chi(d):=\left(\frac{(-1)^k\cdot s}{d}\right)$, where
$s:=\prod_{\delta | N}\delta^{r_{\delta}}$.
Moreover, if $f(z)$ is holomorphic (resp. vanishes) at all of the cusps of $\Gamma_0(N)$, then $f(z)\in M_k(\Gamma_0(N),\chi)$ (resp. $S_k(\Gamma_0(N),\chi)$).

After theorem 2.1, we can easily get the following lemma about $f_b(z)$:

\textbf{Lemma 2.2}  Let $f_b(12z):=\eta(12z)^2\eta(12bz)^2$ as above. Then $f_b(12z)\in S_2(\Gamma_0(144b))$.

We still need the theory of CM forms. The basic idea of CM form can be found in Ribet[5]. Let K be an imaginary quadratic field of discriminant D, and $\varepsilon_K$ be the quadratic character associated to K. Next, let $c$ be a Hecke character of $K$ with exponent $k-1$ and conductor $f_c$. That means, view $c$ as a homomorphism:
\begin{displaymath}
(fractional~ideals~of~K~prime~to~f_c) \rightarrow C^{\ast}
\end{displaymath}
we have
\begin{displaymath}
c((\alpha))=\alpha^{k-1}
\end{displaymath}
for all $\alpha\in F^{\ast}$ such that $\alpha \equiv 1~(mod^{\ast}~f_c)$.

Next, associated to $c$, define a Dirichlet character $\omega_c$ mod $N(f_c)$ as:
\begin{displaymath}
\omega_c(n)=c((n))/n^{k-1}
\end{displaymath}
for all $n\in \mathbb{Z}$ coprime to $f_c$. $N(f_c)$ denotes the norm of $f_c$.

For any $\delta\in \mathbb{Z}^+$, define series $\varphi_{k,c,\delta}$ as:
\begin{displaymath}
\varphi_{K,c,\delta}=\sum_{a}c(a)q^{\delta\cdot N(a)}
\end{displaymath}
where the sum runs through all integral ideals of $K$ coprime to $f_c$. Then, by theorem 3.4 in [5], $\varphi_{K,c,\delta}$ is a cusp form of weight $k$ and
character $\omega_c\cdot \varepsilon_K$ on $\Gamma_0(\delta\cdot|D|\cdot N(f_c))$, moreover, $\varphi_{K,c,\delta}$ is an eigenform for all Hecke operators
$T_p$ in which $p\nmid \delta\cdot|D|\cdot N(f_c)$. So, in order to make $\varphi_{K,c,\delta}$ to be an element in $S_k(\Gamma_0(N),\varepsilon)$, the following two
conditions are necessary and sufficient:
\begin{displaymath}
\delta\cdot|D|\cdot N(f_c)~|~N
\end{displaymath}
\begin{displaymath}
\omega_c\cdot \varepsilon_K=\varepsilon
\end{displaymath}
All possible $\varphi_{K,c,\delta}$ satisfying the above two conditions generate a subspace of $S_k(\Gamma_0(N),\varepsilon)$, which denoted as $S_k^{cm}(\Gamma_0(N),\varepsilon)$.
Then, by Serre[6], an element of $S_k(\Gamma_0(N),\varepsilon)$ is lacunary if and only if it is an element of $S_k^{cm}(\Gamma_0(N),\varepsilon)$.
~\\

\section{Proof of Theorem 1}
~\\

With the above theory, we know that $f_b(12z)$ is lacunary if and only if $f_b(12z)=\sum \varphi_i \in S_2^{cm}(\Gamma_0(144b),\varepsilon)$, where $\varphi_i$ are some suitable
$\varphi_{K,c,\delta}$. First we have the following lemmas:

\textbf{Lemma 3.1} $f_b(12z)=\sum \varphi_i$ as above. $b$ as in Theorem 1. Then as a CM form, $\varphi_i$ can only be associated to one of the following four fields: $\mathbb{Q}(\sqrt{-1})$,
$\mathbb{Q}(\sqrt{-2})$, $\mathbb{Q}(\sqrt{-3})$ and $\mathbb{Q}(\sqrt{-6})$.

\textbf{Proof} Assume $\varphi_i$ is generated by $\mathbb{Q}(\sqrt{D})$, $k$, $f_c$ and $\varepsilon_K$, $\omega_c$. So we have $|D|\cdot N(f_c)~|144b$, and $\omega_c\cdot \varepsilon_K=1_N$,
where $1_N$ is the trivial character on $\Gamma_0(144b)$. By $\omega_c\cdot \varepsilon_K=1_N$, if there is some prime $p\neq 2,3$ such that $p\mid|D|$, that is $p$ is a divisor of the module of character $\varepsilon_K$, so p must be a divisor of the module of character $\omega_c$, too. But we already know that $\omega_c$ is a Dirichlet character mod $N(f_c)$, so there must be
$p\mid N(f_c)$, and $p^2\mid |D|\cdot N(f_c)~|144b$. This is a contradiction to our assumption that b is a square-free integer. This completes the proof. $\Box$

\textbf{Lemma 3.2} $\varphi_i$ is a CM form associated to $\mathbb{Q}(\sqrt{-1})$, $\mathbb{Q}(\sqrt{-2})$, $\mathbb{Q}(\sqrt{-3})$ or $\mathbb{Q}(\sqrt{-6})$. $T_{23}$ is the usual Hecke operator $T_p$ with $p=23$. Then $\varphi_i|T_{23}=0$.

This lemma is obvious because 23 remains prime in each of those four fields, so the 23-rd coefficient of $\varphi_i$ must be zero. Combined with the fact that $\varphi_i$ is an eigenform of
$T_{23}$, and the computing formula of $T_p$: $a_n(T_p(f))=a_{np}(f)+\chi(p)p^{k-1}a_{n/p}(f)$, this lemma can be easily obtained. $\Box$

\textbf{Lemma 3.3}  $f_b(12z)=\sum \varphi_i$ as in Theorem 1. Then if $b\geq 176$, $f_b(12z)$ cannot be lacunary.

\textbf{Proof} By Lemma 3.2, if $f_b(12z)$ is lacunary, then $f_b(12z)|T_{23}=0$. We will show that if $b$ is large enough, $f_b(12z)|T_{23}$ cannot be zero. Denote
$f_b(12z)=q^{1+b}\prod_{n=1}^{\infty}(1-q^{12n})^2(1-q^{12bn})^2=\sum_{n=1}^{\infty}a(n)q^n$ and $\prod_{n=1}^{\infty}(1-q^{n})^2=\sum_{n=1}^{\infty}b(n)q^n$. If $n<12b$, then $a(n)$ has nothing to do with $(1-q^{12bn})^2$, so we have $a(n)=b(\frac{n-(1+b)}{12})$. Denote $T_{23}(f)=\sum_{n=1}^{\infty}c(n)q^n$, Using the formula of $T_p$, we have $c(n)=a(23n)+p\cdot a(n/23)$. When $23\nmid n$, $c(n)=a(23n)=b(\frac{23n-(1+b)}{12})$. In order to make $\frac{23n-(1+b)}{12}$ an integer, we need $n\equiv -(1+b)~(mod~12)$.
If there exists some $c(n)\neq 0$, then $f_b(12z)|T_{23}\neq 0$. In summary, after three conditions: (1) $1+b\leq n<12b$; (2) $23\nmid n$; (3) $n\equiv -(1+b)~(mod~12)$; If there exists such $n$ that $b(\frac{23n-(1+b)}{12})\neq 0$, then $f_b(12z)$ cannot be lacunary.

Let $n_0$ be the smallest $n$ that makes $\frac{23n-(1+b)}{12}$ an positive integer. So $\frac{23n-(1+b)}{12}\in \{1,2,3,\ldots,23\}$, and when $n_0\rightarrow n_0+12$,
$\frac{23n-(1+b)}{12}\rightarrow \frac{23n-(1+b)}{12}+23$. Obviously, in $n_0+12t$, $t=0,1,2,\ldots,22$, these 23 numbers, only one can be divided by 23. We list the first 1000 coefficients of $b(n)$ in appendix 1.

From the table of $b(n)$ above, we only need $t\geq 5$, then there will be at least two numbers not equal to zero in each column, except the 21st column, which we will discuss later. So there must be at least one number that suits $23\nmid n$. In order to suit $1+b\leq n<12b$, we need $b>\frac{12*7*23}{11}$, that is $b\geq 176$. As to the 21st column, in this
condition, $\frac{23n_0-(1+b)}{12}=21$, that is $21n_0-(1+b)=21*12$. Taking module of 23 on the both sides, we get $k\equiv 0~(mod~23)$, which contradicts our assumption in Theorem 1. In summary, we can say that when $b\geq 176$, $f_b(12z)$ cannot be lacunary. $\Box$

Now, in order to prove Theorem 1, we only need to check the lacunarity of each case of $b<176$ by computer. The method is to check whether $f_b(12z)|T_{23}$ (or $T_{47}$, when $23|b$) is zero. At last we find that $f_b(12z)$ cannot be lacunary except $b=1,2,3,4$ or $16$. We will prove those five cases in Proposition 4.1. Theorem 1 is proved.
~\\

\section{Computation when $b=1,2,3,4,16$}
~\\

When b=1,2,3,4,16, $f_b(z)$ is lacunary. So it can be expressed as a linear combination of CM forms. We list these specific expressions in Proposition 4.1:

\textbf{Proposition 4.1}

(1) $b=1$; $f_b(6z)=\eta(6z)^4=\varphi_{K,c}(z)\in S_2(\Gamma_0(36))$, where $K=\mathbb{Q}(\sqrt{-3})$, $f_c=(2\sqrt{-3})$. This can be found in [1].

(2) $b=2$; $f_b(4z)=\eta(4z)^2\eta(8z)^2=\varphi_{K,c}(z)\in S_2(\Gamma_0(32))$, where $K=\mathbb{Q}(\sqrt{-1})$, $f_c=(2(1+i))$.

(3) $b=3$; $f_b(3z)=\eta(3z)^2\eta(9z)^2=\varphi_{K,c}(z)\in S_2(\Gamma_0(27))$, where $K=\mathbb{Q}(\sqrt{-3})$, $f_c=(3)$.

(4) $b=4$; $f_b(12z)=\eta(12z)^2\eta(48z)^2=\frac{1}{8}(\varphi_{K,c_-}(z)-\varphi_{K,c_+}(z))\in S_2(\Gamma_0(576))$, where $K=\mathbb{Q}(\sqrt{-1})$, $f_c=(12)$.

(5) $b=16$; $f_b(12z)=\eta(12z)^2\eta(16\cdot12z)^2=\frac{1}{16}(\varphi_{603}-\varphi_{203})-\frac{\sqrt{1-2\sqrt{-6}}}{16\cdot(6-2\sqrt{-6})}(\varphi_{130}-\varphi_{130}'+\varphi_{310}-\varphi_{310}')
\in S_2(\Gamma_0(144\cdot 16))$, where $\varphi_{603}$ and $\varphi_{203}$ are associated to $K=\mathbb{Q}(\sqrt{-1})$, $f_c=(24)$, and $\varphi_{130}$, $\varphi_{130}'$, $\varphi_{310}$, $\varphi_{310}'$ are associated to $K=\mathbb{Q}(\sqrt{-6})$, $f_c=(4\sqrt{-6})$.

\textbf{Proof}  First we need to compute the specific Hecke character on each case. On case (2), the hecke character $c$ is defined as follows:
choose $a+bi$ to be primary, that is, if $b\equiv 0~(mod~4)$, then $a\equiv 1~(mod~4)$; if $b\equiv 2~(mod~4)$, then $a\equiv 3~(mod~4)$. Then $c((a+bi))=a+bi$. On case (3),
the hecke character $c$ is defined as follows: choose $a+b\omega$ to be primary, that is, $a\equiv 2~(mod~3)$, $b\equiv 0~(mod~3)$. Then $c((a+b\omega))=-(a+b\omega)$. On case (4),
the hecke character $c_+$ and $c_-$ is defined as follows: choose $a+bi$ to be primary, $a+bi\equiv (1-i)^u~(mod~3)$, $u=0,1,\ldots,7$, and $a+bi\equiv (-1+2i)^v~(mod~4)$, $v=0,1$. Then $c_+((a+bi))=i^u\cdot(-1)^v\cdot(a+bi)$, $c_-((a+bi))=(-i)^u\cdot(-1)^v\cdot(a+bi)$. Case (5) is a little bit complicated. To $\varphi_{603}$ and $\varphi_{203}$, the hecke character $c_{603}$ and $c_{203}$ are defined as: choose $a+bi$ to be "standard", which means $a+bi\equiv 1,2+i,3+4i,2+3i,4+i,-1-2i,3i,or~-3-2i~(mod~8)$, if $a+bi\equiv (1-i)^u~(mod~3)$, $u=0,1,\ldots,7$, and $a+bi\equiv (2+i)^v(4+i)^w~(mod~8)$, $v=0,1,2,3$, $w=0,1,2,3$. Then $c_{rst}(a+bi)=\zeta_8^{ru}\cdot i^{sv} \cdot i^{tw}\cdot(a+bi)$, $\zeta_8=e^{\frac{2\pi i}{8}}$. To $\varphi_{130}$, $\varphi_{130}'$, $\varphi_{310}$, $\varphi_{310}'$, the hecke character $c_{130}$, $c_{130}'$, $c_{310}$ and $c_{310}'$ are defined as:
choose $a+b\sqrt{-6}$ to be "standard", which means $a+b\sqrt{-6}\equiv 1,1+\sqrt{-6},1-\sqrt{-6},5,or~the~product~of~them~(mod~4\sqrt{-6})$,
if $a+b\sqrt{-6}\equiv (1+\sqrt{-6})^u\cdot(1-\sqrt{-6})^v\cdot 5^w$, $u,v=0,1,2,3$, $w=0,1$. Then $c_{rst}(a+b\sqrt{-6})=i^{ru}\cdot i^{sv}\cdot (-1)^{tw}\cdot (a+b\sqrt{-6})$. However,
$\mathbb{Q}(\sqrt{-6})$ is not a Principal Ideal Domain, that means not every ideal in $\mathbb{Q}(\sqrt{-6})$ is principal. Let $5=(5,2+\sqrt{-6})\cdot(5,2-\sqrt{-6})$, and denote
$(5,2+\sqrt{-6})$ as $\alpha$. Every non-principal ideal plus $\alpha$ should be a principal ideal because $\mathbb{Q}(\sqrt{-6})$ has Class Number 2. Let $c(\alpha)=\sqrt{c(\alpha^2)}=\sqrt{c(-1+2\sqrt{-6})}$, and $c'(\alpha)=-\sqrt{c(-1+2\sqrt{-6})}$.
That would make the definition explicit.

Next we will verify the above result by Sturm's Theorem[7]:

\textbf{Sturm's Theorem}    Suppose that N is a positive integer, $p$ is prime, and $f(z)$, $g(z)\in M_k(\Gamma_0(N),\chi)\cap \mathbb{Z}[[q]]$. If
\begin{displaymath}
ord_p(f(z)-g(z))>\frac{k}{12}[SL_2(\mathbb{Z}):\Gamma_0(N)]
\end{displaymath}
in which
\begin{displaymath}
[SL_2(\mathbb{Z}):\Gamma_0(N)]=N \prod_{\ell~prime:\ell |N}(1+\frac{1}{\ell})
\end{displaymath}
then $f(z)\equiv g(z)~(mod~p)$.

Case (1)(2)(3) is trivial, because they are in one-dimension linear spaces. To case (4) and (5), although it is very complicated to get the above result, it is relatively easy to check it. By Sturm's Theorem, we only need to verify the first 192 and 768 coefficients, which are listed in appendix 2 and appendix 3. $\Box$

Here we make an example to show how to compute the coefficients of the case $b=16$. Let
$n=29645=5*7^2*11^2$. Denote $\varphi_{uvw}=\sum a_{uvw}(n)q^n$, and $\varphi_{uvw}'=\sum a_{uvw}'(n)q^n$. $5=(2+i)(-1-2i)$ in $\mathbb{Q}(\sqrt{-1})$, so $a_{603}(5)=2$, $a_{203}(5)=-2$.
By the formula $a(p^r)=a(p)a(p^{r-1})-p\cdot a(p^{r-2})$, we can get $a_{603}(49)=a_{203}(49)=-7$, $a_{603}(121)=a_{203}(121)=-11$, so $a_{603}(n)=2*7*11$, $a_{203}(n)=-2*7*11$.
Next, $5=(5,2+\sqrt{-6})\cdot(5,2-\sqrt{-6})$ in $\mathbb{Q}(\sqrt{-6})$, so $a_{130}(5)=a_{310}(5)=\frac{6-2\sqrt{-6}}{\sqrt{1-2\sqrt{-6}}}$, $a_{130}'(5)=a_{310}'(5)=-\frac{6-2\sqrt{-6}}{\sqrt{1-2\sqrt{-6}}}$. There are three
ideals whose norm are 49: $(-5+2\sqrt{-6})$, $(-5-2\sqrt{-6})$ and (7), so $a_{130}(49)=a_{130}'(49)=a_{310}(49)=a_{310}'(49)=17$. There are also three ideals whose norm are 121: $(5+4\sqrt{-6})$, $(5-4\sqrt{-6})$ and (11), so $a_{130}(121)=a_{130}'(121)=a_{310}(121)=a_{310}'(121)=21$. So $a_{130}(29645)=a_{310}(29645)=17*21*\frac{6-2\sqrt{-6}}{\sqrt{1-2\sqrt{-6}}}$, $a_{130}'(29645)=a_{310}'(29645)=-17*21*\frac{6-2\sqrt{-6}}{\sqrt{1-2\sqrt{-6}}}$. Put all these data in the formula, we get $a(29645)=\frac{1}{4}(7*11-17*21)=-70$.
Verifying it by directly computing with computer, we can see this is the right result.
~\\

~\\

\textbf{Acknowledgement}
~\\

The author is grateful to my tutor Chunlei Liu for the guide throughout the article and the help from Haobo Dai.
~\\

\centerline{References}
~\\

[1] J.P.Serre, Sur la lacunarite des puissances de eta, Glasgow Math. J.27(1985). 203-221

[2] E.Clader, Y.Kemper, M.Wage, Lacunarity of certain partition theoretic generating functions, Proceedings of the American Mathematic Society. Volume 137(9),2009. 2959-2968

[3] Y.Martin, Multiplicative eta-quotients, Trans. Amer. Math. Soc. 348(1996). 4825-4856

[4] K.Ono, The Web of Modularity: Arithmetic of the Coefficients of Modular Forms and q-series, CMBS Regional Conference Series in Mathematics, American Mathematical Society,
Providence, RI, 2004

[5] K.A.Ribet, Galois representations attached to eigenforms with Nebentypus

[6] J.P.Serre, Quelques applications du theoreme de densite de Chebatorev, Publ. Math. I.H.E.S.54 (1981). 123-201

[7] John J. Webb, Arithmetic of the 13-regular partition function modulo 3, Ramanujan J(2011)25:49-56

\newpage

\centerline{\textbf{Appendix 1}}

This is the first 1000 coefficients of $\prod_{n=1}^{\infty}(1-q^{n})^2=\sum_{n=1}^{\infty}b(n)q^n$. The first row, from left to right, are $b(1)$, $b(2)$, $\ldots$, $b(23)$, the second row are  $b(24)$, $b(25)$, $\ldots$, $b(46)$, and so on. We can see that in the first six rows, except for the 21st column, there are at least two numbers not equal to zero in each column. That's just what we need in Lemma 3.3.

\begin{displaymath}
\begin{tabular}{*{23}{c}}
-2&-1&2&1&2&-2&0&-2&-2&1&0&0&2&3&-2&2&0&0&-2&-2&0&0&-2\\
-1&0&2&2&-2&2&1&2&0&2&-2&-2&2&0&-2&0&-4&0&0&0&1&-2&0\\
0&2&0&2&2&1&-2&0&2&2&0&0&-2&0&-2&0&-2&2&0&-4&0&0&-2\\
-1&2&0&2&0&0&0&-2&2&4&1&0&0&2&-2&2&-2&0&0&2&0&-2&0\\
-2&-2&0&-2&0&0&0&2&-2&-1&-2&-2&4&0&0&2&0&2&0&0&0&3&-2\\
0&0&4&2&0&-2&0&0&-2&0&0&-2&0&2&0&-2&0&-2&-2&-2&0&0&2\\
-2&-1&2&0&0&0&2&-2&0&-2&2&2&2&2&0&1&2&-2&0&0&0&0&2\\
0&0&0&0&2&0&-2&-2&-4&2&0&0&-2&0&-2&0&2&0&-2&0&0&-4&1\\
-2&0&0&4&2&-2&2&0&0&0&-2&4&0&-2&2&1&0&2&2&0&0&2&0\\
0&-4&2&0&0&-2&0&0&2&0&-2&0&0&0&0&-2&-2&-4&-2&2&0&2&0\\
0&0&-2&-1&0&2&0&-2&0&0&0&0&2&0&0&2&0&4&2&-2&0&-1&2\\
-2&2&0&0&0&2&-2&4&0&0&2&-2&0&0&-2&-2&0&-2&0&0&0&2&-2\\
0&0&0&-2&-2&0&0&0&0&-2&0&-2&-2&1&0&0&2&2&2&0&0&-2&-4\\
4&2&0&-2&0&0&2&0&2&2&3&-2&0&2&2&0&-2&0&0&0&0&-2&0\\
2&2&0&0&-2&0&0&0&0&0&0&0&-2&-4&0&2&-4&0&-2&0&0&2&0\\
2&0&-2&0&-2&0&-3&0&0&2&-2&0&2&0&0&0&0&2&0&4&0&0&0\\
0&2&-4&0&2&1&0&2&2&-4&2&2&0&0&2&0&2&0&0&-2&0&0&-2\\
0&0&-2&0&-2&0&0&2&2&2&-2&0&-4&-2&0&0&0&-2&0&-2&0&-2&2\\
0&-2&0&0&0&1&0&0&2&0&-2&2&0&0&0&-4&0&0&2&2&0&2&0\\
0&0&2&0&0&4&3&-2&0&0&0&0&0&-2&0&-2&2&0&4&0&0&0&2\\
0&0&-2&-2&2&0&0&0&-4&-2&2&0&0&-2&-2&0&2&-2&-2&0&0&0&0\\
4&-2&0&0&-2&-2&-2&0&0&-4&1&4&-2&0&0&0&0&0&2&2&0&0&2\\
0&0&0&-2&2&0&0&0&0&0&0&-2&2&1&6&0&2&0&-2&0&0&2&0\\
-2&2&0&2&0&-2&0&0&0&-2&-2&0&0&0&2&0&2&0&-2&0&0&0&0\\
2&-4&-2&-2&0&-4&2&0&-2&0&0&0&-2&2&0&0&2&-2&0&0&-2&1&0\\
0&2&0&-2&-2&0&0&-2&0&0&4&0&2&-2&4&0&0&0&-2&4&0&0&2\\
0&0&2&1&-2&0&0&0&0&-2&2&2&-2&0&0&2&2&2&0&2&0&0&0\\
0&-2&-2&-4&0&0&2&2&-2&-2&0&0&0&-2&0&0&0&-2&2&2&0&-2&0\\
-2&-2&0&0&2&0&-4&-2&0&0&0&-2&0&0&0&-1&2&0&4&-4&0&2&0\\
2&0&0&-2&0&0&2&2&0&-2&0&0&-2&0&0&0&2&2&2&0&0&2&3\\
2&0&2&-2&0&-2&0&-2&2&0&0&2&0&0&2&-4&0&0&0&2&0&0&0\\
0&2&-4&0&0&2&2&0&0&-2&-2&0&-2&-2&2&-4&-2&0&0&0&0&2&-2\\
0&0&-2&0&2&0&-4&2&2&0&0&0&0&-2&-2&-1&0&-2&0&0&0&0&2\\
0&-2&0&2&0&2&0&0&4&0&2&0&0&-2&0&0&0&-2&2&2&0&0&0\\
2&3&-2&-2&0&2&0&0&0&0&0&4&0&0&0&4&0&2&-2&0&0&-2&-2\\
0&0&-2&0&-2&2&0&0&-2&2&0&0&-2&-2&4&0&0&0&0&-2&0&-2&0\\
2&0&-2&0&-4&-2&0&0&0&0&-2&0&-2&0&0&0&2&0&0&2&0&-1&0\\
2&0&-4&0&-2&0&0&6&2&-2&0&-2&-2&0&0&0&0&2&0&2&0&-2&0\\
\end{tabular}
\end{displaymath}
\begin{displaymath}
\begin{tabular}{*{23}{c}}
2&2&2&0&0&0&4&0&2&1&0&0&0&-2&0&0&-2&0&0&-2&0&2&0\\
4&2&0&0&-2&2&0&-2&0&0&0&2&0&2&0&-2&0&0&-2&0&0&0&-4\\
0&0&0&-4&0&2&0&0&-2&-2&2&0&2&0&-2&0&0&-2&0&0&0&0&-2\\
2&0&-2&0&0&0&-4&0&2&-2&1&-4&0&0&2&0&0&0&0&0&0&4&2\\
0&2&0&0&0&0&-2&2&-2&0&0&-2&0&0&4&0&0&-2&2&2&0&0&0\\
3&0&0&0&0&2&2&0&0&0&0&
\end{tabular}
\end{displaymath}

\newpage

\centerline{\textbf{Appendix 2}}

Let $f_4(12z)=\sum a(n)q^n$, and $\varphi_{K,c_+}(z)=\sum b(n)q^n$, $\varphi_{K,c_-}(z)=\sum c(n)q^n$. By proposition 4.1(4), we only need to verify $a(n)=\frac{1}{8}(c(n)-b(n))$.

We omit those $n$ that $2|n$ or $3|n$, because those coefficients are all zero.
\begin{displaymath}
\begin{array}{cccc}
\quad n\quad&\quad a(n)\quad&\quad b(n)\quad&\quad c(n)\quad\\
1&0&1&1\\
5&1&-4&4\\
7&0&0&0\\
11&0&0&0\\
13&0&6&6\\
17&-2&8&-8\\
19&0&0&0\\
23&0&0&0\\
25&0&11&11\\
29&-1&4&-4\\
31&0&0&0\\
35&0&0&0\\
37&0&2&2\\
41&2&-8&8\\
43&0&0&0\\
47&0&0&0\\
49&0&-7&-7\\
53&-1&4&-4\\
55&0&0&0\\
59&0&0&0\\
61&0&10&10\\
65&6&-24&24\\
67&0&0&0\\
71&0&0&0\\
73&0&6&6\\
77&0&0&0\\
79&0&0&0\\
83&0&0&0\\
85&0&-32&-32\\
89&-4&16&-16\\
91&0&0&0\\
95&0&0&0\\
97&0&-18&-18\\
101&-5&20&-20\\
103&0&0&0\\
107&0&0&0\\
109&0&-18&-18\\
113&-4&16&-16\\
\end{array}
\end{displaymath}
\begin{displaymath}
\begin{array}{cccc}
\quad n\quad&\quad a(n)\quad&\quad b(n)\quad&\quad c(n)\quad\\
115&0&0&0\\
119&0&0&0\\
121&0&-11&-11\\
125&6&-24&24\\
127&0&0&0\\
131&0&0&0\\
133&0&0&0\\
137&-2&8&-8\\
139&0&0&0\\
143&0&0&0\\
145&0&-16&-16\\
149&5&-20&20\\
151&0&0&0\\
155&0&0&0\\
157&0&-22&-22\\
161&0&0&0\\
163&0&0&0\\
167&0&0&0\\
169&0&23&23\\
173&1&-4&4\\
175&0&0&0\\
179&0&0&0\\
181&0&-18&-18\\
185&2&-8&8\\
187&0&0&0\\
191&0&0&0\\
\end{array}
\end{displaymath}

\newpage

\centerline{\textbf{Appendix 3}}

Let $f_{16}(12z)=\sum a(n)q^n$, and $\varphi_{603}(z)=\sum b_1(n)q^n$, $\varphi_{203}(z)=\sum b_2(n)q^n$, $\frac{\sqrt{1-2\sqrt{-6}}}{6-2\sqrt{-6}}\varphi_{130}(z)=\sum c_1(n)q^n$,
$\frac{\sqrt{1-2\sqrt{-6}}}{6-2\sqrt{-6}}\varphi_{130}'(z)=\sum c_2(n)q^n$, $\frac{\sqrt{1-2\sqrt{-6}}}{6-2\sqrt{-6}}\varphi_{310}(z)=\sum c_3(n)q^n$,
$\frac{\sqrt{1-2\sqrt{-6}}}{6-2\sqrt{-6}}\varphi_{310}'(z)=\sum c_4(n)q^n$.
Denote $\frac{\sqrt{1-2\sqrt{-6}}}{6-2\sqrt{-6}}=t$. By proposition 4.1(5), we need to verify that $a(n)=\frac{1}{16}\left(b_1(n)-b_2(n)\right)-\frac{1}{16}\left(c_1(n)-c_2(n)+c_3(n)-c_4(n)\right)$.

We omit those $n$ that $2|n$ or $3|n$, because those coefficients are all zero.

"$\ast$" means that because of the "minus" in the formula, "$*-*=0$", so we don't need to compute these specific values. These case are:

(1) If $p=(a+bi)$ is an ideal in $\mathbb{Q}(\sqrt{-1})$, $3|a$ or $3|b$, then $c_+(p)=c_-(p)$.

(2) If $p$ is a principle ideal in $\mathbb{Q}(\sqrt{-6})$, then $c_{130}(p)=c_{130}'(p)$, $c_{310}(p)=c_{310}'(p)$.

\begin{displaymath}
\begin{array}{cccccccc}
\quad n \quad&\quad a(n)\quad & \quad b_1(n) \quad & \quad b_2(n)\quad & \quad c_1(n)\quad & \quad c_2(n) \quad & \quad c_3(n)\quad & \quad c_4(n)\quad\\
1&0&1&1&1&1&1&1\\
5&0&2&-2&1&-1&1&-1\\
7&0&0&0&-2\sqrt{6}t&-2\sqrt{6}t&2\sqrt{6}t&2\sqrt{6}t\\
11&0&0&0&\frac{4}{\sqrt{6}}&-\frac{4}{\sqrt{6}}&-\frac{4}{\sqrt{6}}&\frac{4}{\sqrt{6}}\\
13&0&4&4&0&0&0&0\\
17&1&8&-8&0&0&0&0\\
19&0&0&0&0&0&0&0\\
23&0&0&0&0&0&0&0\\
25&0&-1&-1&7t&7t&7t&7t\\
29&-2&-10&10&3&-3&3&-3\\
31&0&0&0&-2\sqrt{6}t&-2\sqrt{6}t&2\sqrt{6}t&2\sqrt{6}t\\
35&0&0&0&-2\sqrt{6}&2\sqrt{6}&2\sqrt{6}&-2\sqrt{6}\\
37&0&12&12&0&0&0&0\\
41&-1&-8&8&0&0&0&0\\
43&0&0&0&0&0&0&0\\
47&0&0&0&0&0&0&0\\
49&0&-7&-7&17t&17t&17t&17t\\
53&2&14&-14&-1&1&-1&1\\
55&0&0&0&\ast&\ast&\ast&\ast\\
59&0&0&0&\frac{8}{\sqrt{6}}&-\frac{8}{\sqrt{6}}&-\frac{8}{\sqrt{6}}&\frac{8}{\sqrt{6}}\\
61&0&12&12&0&0&0&0\\
65&1&8&-8&0&0&0&0\\
67&0&0&0&0&0&0&0\\
71&0&0&0&0&0&0&0\\
73&0&6&6&14t&14t&14t&14t\\
77&2&0&0&-8&8&-8&8\\
79&0&0&0&6\sqrt{6}t&6\sqrt{6}t&-6\sqrt{6}t&-6\sqrt{6}t\\
83&0&0&0&-\frac{4}{\sqrt{6}}&\frac{4}{\sqrt{6}}&\frac{4}{\sqrt{6}}&-\frac{4}{\sqrt{6}}\\
85&0&16&16&0&0&0&0\\
\end{array}
\end{displaymath}
\begin{displaymath}
\begin{array}{cccccccc}
\quad n \quad&\quad a(n)\quad & \quad b_1(n) \quad & \quad b_2(n)\quad & \quad c_1(n)\quad & \quad c_2(n) \quad & \quad c_3(n)\quad & \quad c_4(n)\quad\\
89&-2&-16&16&0&0&0&0\\
91&0&0&0&0&0&0&0\\
95&0&0&0&0&0&0&0\\
97&0&18&18&2t&2t&2t&2t\\
101&0&-2&2&-1&1&-1&1\\
103&0&0&0&-6\sqrt{6}t&-6\sqrt{6}t&6\sqrt{6}t&6\sqrt{6}t\\
107&0&0&0&\frac{8}{\sqrt{6}}&-\frac{8}{\sqrt{6}}&-\frac{8}{\sqrt{6}}&\frac{8}{\sqrt{6}}\\
109&0&20&20&0&0&0&0\\
113&-2&-16&16&0&0&0&0\\
115&0&0&0&0&0&0&0\\
119&0&0&0&0&0&0&0\\
121&0&-11&-11&21t&21t&21t&21t\\
125&-2&-12&12&2&-2&2&-2\\
127&0&0&0&2\sqrt{6}t&2\sqrt{6}t&-2\sqrt{6}t&-2\sqrt{6}t\\
131&0&0&0&\frac{16}{\sqrt{6}}&-\frac{16}{\sqrt{6}}&-\frac{16}{\sqrt{6}}&\frac{16}{\sqrt{6}}\\
133&0&0&0&0&0&0&0\\
137&1&8&-8&0&0&0&0\\
139&0&0&0&0&0&0&0\\
143&0&0&0&0&0&0&0\\
145&0&-20&-20&36t&36t&36t&36t\\
149&0&-14&14&-7&7&-7&7\\
151&0&0&0&-10\sqrt{6}t&-10\sqrt{6}t&10\sqrt{6}t&10\sqrt{6}t\\
155&0&0&0&-2\sqrt{6}&2\sqrt{6}&2\sqrt{6}&-2\sqrt{6}\\
157&0&12&12&0&0&0&0\\
161&0&0&0&0&0&0&0\\
163&0&0&0&0&0&0&0\\
167&0&0&0&0&0&0&0\\
169&0&3&3&\ast&\ast&\ast&\ast\\
173&2&26&-26&5&-5&5&-5\\
175&0&0&0&\ast&\ast&\ast&\ast\\
179&0&0&0&-\frac{8}{\sqrt{6}}&\frac{8}{\sqrt{6}}&\frac{8}{\sqrt{6}}&-\frac{8}{\sqrt{6}}\\
181&0&20&20&0&0&0&0\\
185&3&24&-24&0&0&0&0\\
187&0&0&0&0&0&0&0\\
191&0&0&0&0&0&0&0\\
193&0&14&14&\ast&\ast&\ast&\ast\\
197&-2&-2&2&7&-7&7&-7\\
199&0&0&0&\ast&\ast&\ast&\ast\\
203&0&0&0&-6\sqrt{6}&6\sqrt{6}&6\sqrt{6}&-6\sqrt{6}\\
205&0&-16&-16&0&0&0&0\\
209&0&0&0&0&0&0&0\\
211&0&0&0&0&0&0&0\\
215&0&0&0&0&0&0&0\\
217&0&0&0&\ast&\ast&\ast&\ast\\
221&4&32&-32&0&0&0&0\\
223&0&0&0&\ast&\ast&\ast&\ast\\
\end{array}
\end{displaymath}
\begin{displaymath}
\begin{array}{cccccccc}
\quad n \quad&\quad a(n)\quad & \quad b_1(n) \quad & \quad b_2(n)\quad & \quad c_1(n)\quad & \quad c_2(n) \quad & \quad c_3(n)\quad & \quad c_4(n)\quad\\
227&0&0&0&\frac{20}{\sqrt{6}}&-\frac{20}{\sqrt{6}}&-\frac{20}{\sqrt{6}}&\frac{20}{\sqrt{6}}\\
229&0&4&4&0&0&0&0\\
233&2&16&-16&0&0&0&0\\
235&0&0&0&0&0&0&0\\
239&0&0&0&0&0&0&0\\
241&0&-30&-30&\ast&\ast&\ast&\ast\\
245&-6&-14&14&17&-17&17&-17\\
247&0&0&0&0&0&0&0\\
251&0&0&0&\frac{4}{\sqrt{6}}&-\frac{4}{\sqrt{6}}&-\frac{4}{\sqrt{6}}&\frac{4}{\sqrt{6}}\\
253&0&0&0&0&0&0&0\\
257&-4&-32&32&0&0&0&0\\
259&0&0&0&0&0&0&0\\
263&0&0&0&0&0&0&0\\
265&0&-28&-28&\ast&\ast&\ast&\ast\\
269&-4&-26&26&-3&3&-3&3\\
271&0&0&0&\ast&\ast&\ast&\ast\\
275&0&0&0&\frac{28}{\sqrt{6}}&-\frac{28}{\sqrt{6}}&-\frac{28}{\sqrt{6}}&\frac{28}{\sqrt{6}}\\
277&0&\ast&\ast&0&0&0&0\\
281&4&32&-32&0&0&0&0\\
283&0&0&0&0&0&0&0\\
287&0&0&0&0&0&0&0\\
289&0&47&47&-17&-17&-17&-17\\
293&-2&-34&34&-9&9&-9&9\\
295&0&0&0&\ast&\ast&\ast&\ast\\
299&0&0&0&0&0&0&0\\
301&0&0&0&0&0&0&0\\
305&3&24&-24&0&0&0&0\\
307&0&0&0&0&0&0&0\\
311&0&0&0&0&0&0&0\\
313&0&\ast&\ast&\ast&\ast&\ast&\ast\\
317&4&22&-22&-5&5&-5&5\\
319&0&0&0&\ast&\ast&\ast&\ast\\
323&0&0&0&0&0&0&0\\
325&0&-4&-4&0&0&0&0\\
329&0&0&0&0&0&0&0\\
331&0&0&0&0&0&0&0\\
335&0&0&0&0&0&0&0\\
337&0&\ast&\ast&0&0&0&0\\
341&2&0&0&-8&8&-8&8\\
343&0&0&0&\ast&\ast&\ast&\ast\\
347&0&0&0&-\frac{20}{\sqrt{6}}&\frac{20}{\sqrt{6}}&\frac{20}{\sqrt{6}}&-\frac{20}{\sqrt{6}}\\
349&0&\ast&\ast&0&0&0&0\\
353&-2&-16&16&0&0&0&0\\
355&0&0&0&0&0&0&0\\
359&0&0&0&0&0&0&0\\
361&0&-19&-19&\ast&\ast&\ast&\ast\\
365&-2&12&-12&14&-14&14&-14\\
\end{array}
\end{displaymath}
\begin{displaymath}
\begin{array}{cccccccc}
\quad n \quad&\quad a(n)\quad & \quad b_1(n) \quad & \quad b_2(n)\quad & \quad c_1(n)\quad & \quad c_2(n) \quad & \quad c_3(n)\quad & \quad c_4(n)\quad\\
367&0&0&0&\ast&\ast&\ast&\ast\\
371&0&0&0&2\sqrt{6}&-2\sqrt{6}&-2\sqrt{6}&2\sqrt{6}\\
373&0&\ast&\ast&0&0&0&0\\
377&-5&-40&40&0&0&0&0\\
379&0&0&0&0&0&0&0\\
383&0&0&0&0&0&0&0\\
385&0&0&0&\ast&\ast&\ast&\ast\\
389&6&34&-34&-7&7&-7&7\\
391&0&0&0&0&0&0&0\\
395&0&0&0&6\sqrt{6}&-6\sqrt{6}&-6\sqrt{6}&6\sqrt{6}\\
397&0&\ast&\ast&0&0&0&0\\
401&-5&-40&40&0&0&0&0\\
403&0&0&0&0&0&0&0\\
407&0&0&0&0&0&0&0\\
409&0&\ast&\ast&\ast&\ast&\ast&\ast\\
413&4&0&0&-16&16&-16&16\\
415&0&0&0&\ast&\ast&\ast&\ast\\
419&0&0&0&-\frac{28}{\sqrt{6}}&\frac{28}{\sqrt{6}}&\frac{28}{\sqrt{6}}&-\frac{28}{\sqrt{6}}\\
421&0&\ast&\ast&0&0&0&0\\
425&-1&-8&8&0&0&0&0\\
427&0&0&0&0&0&0&0\\
431&0&0&0&0&0&0&0\\
433&0&\ast&\ast&\ast&\ast&\ast&\ast\\
437&0&0&0&0&0&0&0\\
439&0&0&0&\ast&\ast&\ast&\ast\\
443&0&0&0&-\frac{20}{\sqrt{6}}&\frac{20}{\sqrt{6}}&\frac{20}{\sqrt{6}}&-\frac{20}{\sqrt{6}}\\
445&0&-32&-32&0&0&0&0\\
449&5&40&-40&0&0&0&0\\
451&0&0&0&0&0&0&0\\
455&0&0&0&0&0&0&0\\
457&0&\ast&\ast&\ast&\ast&\ast&\ast\\
461&-2&-38&38&-11&11&-11&11\\
463&0&0&0&\ast&\ast&\ast&\ast\\
467&0&0&0&-\frac{28}{\sqrt{6}}&\frac{28}{\sqrt{6}}&\frac{28}{\sqrt{6}}&-\frac{28}{\sqrt{6}}\\
469&0&0&0&0&0&0&0\\
473&0&0&0&0&0&0&0\\
475&0&0&0&0&0&0&0\\
479&0&0&0&0&0&0&0\\
481&0&48&48&0&0&0&0\\
485&4&36&-36&2&-2&2&-2\\
487&0&0&0&\ast&\ast&\ast&\ast\\
491&0&0&0&-\frac{16}{\sqrt{6}}&\frac{16}{\sqrt{6}}&\frac{16}{\sqrt{6}}&-\frac{16}{\sqrt{6}}\\
493&0&-80&-80&0&0&0&0\\
497&0&0&0&0&0&0&0\\
499&0&0&0&0&0&0&0\\
503&0&0&0&0&0&0&0\\
\end{array}
\end{displaymath}
\begin{displaymath}
\begin{array}{cccccccc}
\quad n \quad&\quad a(n)\quad & \quad b_1(n) \quad & \quad b_2(n)\quad & \quad c_1(n)\quad & \quad c_2(n) \quad & \quad c_3(n)\quad & \quad c_4(n)\quad\\
505&0&-4&-4&\ast&\ast&\ast&\ast\\
509&2&-10&10&-13&13&-13&13\\
511&0&0&0&\ast&\ast&\ast&\ast\\
515&0&0&0&-6\sqrt{6}&6\sqrt{6}&6\sqrt{6}&-6\sqrt{6}\\
517&0&0&0&0&0&0&0\\
521&-5&-40&40&0&0&0&0\\
523&0&0&0&0&0&0&0\\
527&0&0&0&0&0&0&0\\
529&0&-23&-23&\ast&\ast&\ast&\ast\\
533&-4&-32&32&0&0&0&0\\
535&0&0&0&\ast&\ast&\ast&\ast\\
539&0&0&0&\frac{68}{\sqrt{6}}&-\frac{68}{\sqrt{6}}&-\frac{68}{\sqrt{6}}&\frac{68}{\sqrt{6}}\\
541&0&\ast&\ast&0&0&0&0\\
545&5&40&-40&0&0&0&0\\
547&0&0&0&0&0&0&0\\
551&0&0&0&0&0&0&0\\
553&0&0&0&\ast&\ast&\ast&\ast\\
557&-8&-38&38&13&-13&13&-13\\
559&0&0&0&0&0&0&0\\
563&0&0&0&\frac{20}{\sqrt{6}}&-\frac{20}{\sqrt{6}}&-\frac{20}{\sqrt{6}}&\frac{20}{\sqrt{6}}\\
565&0&-32&-32&0&0&0&0\\
569&-5&-40&40&0&0&0&0\\
571&0&0&0&0&0&0&0\\
575&0&0&0&0&0&0&0\\
577&0&\ast&\ast&\ast&\ast&\ast&\ast\\
581&-2&0&0&8&-8&8&-8\\
583&0&0&0&\ast&\ast&\ast&\ast\\
587&0&0&0&\frac{32}{\sqrt{6}}&-\frac{32}{\sqrt{6}}&-\frac{32}{\sqrt{6}}&\frac{32}{\sqrt{6}}\\
589&0&0&0&0&0&0&0\\
593&2&16&-16&0&0&0&0\\
595&0&0&0&0&0&0&0\\
599&0&0&0&0&0&0&0\\
601&0&\ast&\ast&\ast&\ast&\ast&\ast\\
605&-8&-22&22&21&-21&21&-21\\
607&0&0&0&\ast&\ast&\ast&\ast\\
611&0&0&0&0&0&0&0\\
613&0&\ast&\ast&0&0&0&0\\
617&4&32&-32&0&0&0&0\\
619&0&0&0&0&0&0&0\\
623&0&0&0&0&0&0&0\\
625&0&-19&-19&\ast&\ast&\ast&\ast\\
629&12&96&-96&0&0&0&0\\
631&0&0&0&\ast&\ast&\ast&\ast\\
635&0&0&0&2\sqrt{6}&-2\sqrt{6}&-2\sqrt{6}&2\sqrt{6}\\
641&1&8&-8&0&0&0&0\\
643&0&0&0&0&0&0&0\\
647&0&0&0&0&0&0&0\\
\end{array}
\end{displaymath}
\begin{displaymath}
\begin{array}{cccccccc}
\quad n \quad&\quad a(n)\quad & \quad b_1(n) \quad & \quad b_2(n)\quad & \quad c_1(n)\quad & \quad c_2(n) \quad & \quad c_3(n)\quad & \quad c_4(n)\quad\\
649&0&0&0&\ast&\ast&\ast&\ast\\
653&2&26&-26&5&-5&5&-5\\
655&0&0&0&\ast&\ast&\ast&\ast\\
659&0&0&0&-\frac{32}{\sqrt{6}}&\frac{32}{\sqrt{6}}&\frac{32}{\sqrt{6}}&-\frac{32}{\sqrt{6}}\\
661&0&\ast&\ast&0&0&0&0\\
665&0&0&0&0&0&0&0\\
667&0&0&0&0&0&0&0\\
671&0&0&0&0&0&0&0\\
673&0&\ast&\ast&\ast&\ast&\ast&\ast\\
677&4&2&-2&-15&15&-15&15\\
679&0&0&0&\ast&\ast&\ast&\ast\\
637&0&-28&-28&0&0&0&0\\
683&0&0&0&-15i&15i&15i&-15i\\
685&0&16&16&0&0&0&0\\
689&7&56&-56&0&0&0&0\\
691&0&0&0&0&0&0&0\\
695&0&0&0&0&0&0&0\\
697&0&-64&-64&0&0&0&0\\
701&-4&-10&10&11&-11&11&-11\\
703&0&0&0&0&0&0&0\\
707&0&0&0&2\sqrt{6}&-2\sqrt{6}&-2\sqrt{6}&2\sqrt{6}\\
709&0&\ast&\ast&0&0&0&0\\
713&0&0&0&0&0&0&0\\
715&0&0&0&0&0&0&0\\
719&0&0&0&0&0&0&0\\
721&0&0&0&\ast&\ast&\ast&\ast\\
725&-4&10&-10&21&-21&21&-21\\
727&0&0&0&\ast&\ast&\ast&\ast\\
731&0&0&0&0&0&0&0\\
733&0&\ast&\ast&0&0&0&0\\
737&0&0&0&0&0&0&0\\
739&0&0&0&0&0&0&0\\
743&0&0&0&0&0&0&0\\
745&0&-28&-28&\ast&\ast&\ast&\ast\\
749&4&0&0&\ast&\ast&\ast&\ast\\
751&0&0&0&\ast&\ast&\ast&\ast\\
755&0&0&0&-10\sqrt{6}&10\sqrt{6}&10\sqrt{6}&-10\sqrt{6}\\
757&0&\ast&\ast&0&0&0&0\\
761&5&40&-40&0&0&0&0\\
763&0&0&0&0&0&0&0\\
767&0&0&0&0&0&0&0\\
\end{array}
\end{displaymath}

\end{document}